\documentclass{owrart}

\usepackage{amsthm,amsmath}
\usepackage{amssymb,tikz}
\usepackage[all]{xy}
\usepackage{graphicx}
\usepackage{url}
\usepackage{mathtools}
\usepackage{hyperref}
\hypersetup{
    colorlinks = true,
    linkcolor = blue,
    urlcolor = blue,
    citecolor = blue,
    anchorcolor = black,
}
\usepackage{color}
\newcommand{\ba}{\begin{align}}
\newcommand{\ea}{\end{align}}
\newcommand{\be}{\begin{equation}}
\newcommand{\ee}{\end{equation}}
\newcommand{\beq}{\begin{eqnarray}}
\newcommand{\eeq}{\end{eqnarray}}
\newcommand{\beqs}{\begin{eqnarray*}}
\newcommand{\eeqs}{\end{eqnarray*}}

\newcommand{\size}[1]{\fontsize{10pt}{\baselineskip}\selectfont{#1}}
\newcommand{\FS}{\mathfrak{F}_{s}}

\begin{document}


\begin{talk}{Arthur Jaffe}
{Reflection Positivity Then and Now}
{Jaffe, Arthur}
\section*{Dedicated to the memory of Robert Schrader}
Konrad Osterwalder and Robert Schrader discovered \textit{reflection positivity} (RP)  in the summer of 1972 at Harvard University.\footnote{This is based on the opening talk on November 20, 2017 at the conference ``Reflection Positivity,'' held at the  Mathematical Research Institute, Oberwolfach, Germany.}
  At the time they were both my postdoctoral fellows, and together we formed the core of a group at Harvard who studied quantum field theory from a mathematical point of view.  
  
RP has since blossomed into an active research area; this conference in Oberwolfach marks the $45^{\rm th}$ anniversary of its discovery.  When Palle Jorgensen, Karl-Hermann Neeb,  Gestur \'Olafsson, and I first envisioned this meeting in 2014, we had hoped  that both Robert and Konrad would be here.  Unfortunately neither is: Robert died in November 2015, and Konrad had to cancel because of a conflict.   
So I have been called upon to make some remarks to set the stage. 

These comments not only relate to the early discovery of RP, but they also  illustrate that  RP is still an active and interesting area of research.   
The original discovery of RP arose from an effort to relate two different mathematical subjects in a specific way:  what property is  needed to start from a classical probability theory of fields (i.e. a statistical mechanics of random fields), and end up with a quantum theory of fields? Of course one wanted to have a framework that would apply to the putative quantum fields that physicists believe describe the interactions of elementary particles.   

The original work focused into finding a set of axioms for Euclidean-covariant Green's functions, that are equivalent to the Wightman axioms for the vacuum expectation values of relativistic quantum fields~\cite{OS1,OS2}.   The Euclidean Green's functions are objects in classical probability theory, while  vacuum expectation values describe the quantum mechanics of fields.  

One can regard  the relation between the  analytically continued expectations (in time) to probability theory,  as generalizing the ``Feynman-Kac formula.'' This formula gives a Wiener integral representation that provides solutions to  the heat-diffusion equation, namely  the analytic continuation of the Schr\"odinger equation.  One side of the equation is purely classical; the other side is an analytic continuation (justified by the positivity of the energy) of quantum theory. 

An immediate application of RP was to provide a framework for the first existence proof for examples of non-linear quantum field theories with scattering~\cite{GJS-Axioms}. These are quantum,  non-linear, scalar wave equations with energy density $\lambda \mathcal{P}(\varphi)_{2}$. Here $\lambda\mathcal{P}$ is a polynomial bounded from below, and the subscript $2$ denotes two-dimensional space time.  \textbf{Theorem:} \textit{For $0<\lambda$ sufficiently small, these models exist and satisfy all the Wightman axioms for a quantum field theory.} The models actually satisfy as well the more-detailed Haag-Ruelle axioms for scattering, which include the existence of an lower and upper mass gap, leading to an isolated eigenvalue $m$ in the mass spectrum. In addition they yield local algebras that satisfy the Haag-Kastler axioms.  For these reasons the discovery of RP marked a turning point in the study of relativistic quantum physics. 
 
Amazingly, the RP property not only arises in quantum physics. But RP connects to areas of research in fields ranging from operator algebras, mathematical analysis, probability, and representation theory on one hand,  to statistical physics and the study of phase transitions on the other.  RP has enabled the proofs of numerous interesting and deep results in far-flung areas of mathematics and physics.  Similar positivity conditions appear in other subjects, so one might dream that new areas of relevance for reflection positivity will appear in the future, in other mathematical areas.

Since RP is now an enormous field, I apologize in advance for my sparse list of references;  these are only meant to be personal and impressionistic, and they cover a small selection of papers  from an enormous universe of possibilities. However, I am confident that the other speakers during this week will fill the gap by citing many other important papers.

 \section{The Background Story}
The analytic continuation of the expectations of fields to imaginary time was proved in the general framework of Wightman in 1956~\cite{Wightman}. This result is a consequence of the assumption that the spectrum of the energy $H$ is non-negative, and that the time-translation automorphism arises from the unitary group $e^{itH}$ acting on the Hilbert space $\mathcal{H}$ of quantum theory.  Furthermore one assumes that the Poincar\'e group acting on $\mathcal{H}$ has an invariant vector (a vacuum). 

This analytic continuation is the non-perturbative parallel to ``Wick rotation'' of time introduced in the physics literature to make sense of the terms in the perturbation expansion of a Lorentz-invariant theory~\cite{Wick}.   Jost studied properties of analytically continued Wightman functions to Schwinger functions, including their symmetry at imaginary time, as detailed in his 1957 book~\cite{Jost}. Schwinger emphasized the notion of Euclidean quantum field theory in his well-known 1958 paper~\cite{SchwingerEQFT}.  According to Miller~\cite{Miller}, Schwinger claimed to have discovered these results about Euclidean fields seven years before he published them. 

In any case, the advent of Euclidean expectation values and Euclidean fields set the stage for the pioneering work of Kurt Symanzik, who proposed in 1964 that the  Euclidean field  could be formulated  as the classical Markov field.  Symanzik described these ideas in preliminary form in his widely-circulated Courant Institute Report~~\cite{Symanzik-NYU} and in the paper \cite{SymanzikJMP}.  He also presented a monumental course in his 1968 Varenna summer school lectures~\cite{Symanzik-Varenna}. 

The Varenna lectures  mark the end of Symanzik's research on Euclidian quantum field theory, which he explains in his introductory lecture. There he estimates the efficacy of the different approaches known at the time for constructing an example of an interacting field. Unfortunately Symanzik under-estimated the value of his own approach. In fact when it was understood some five years later, it turned out to play an important role.

Symanzik's approach cried out for building a  mathematical foundation by which to understand it. His vision appeared to be a promising blueprint for finding a non-trivial field theory.  Symanzik made a valiant attempt to implement  his program for  the $\varphi^{4}$ theory in two space-time dimensions, by working on this problem with S.R.S.~Varadhan.  But the gap between their results (described in the appendix to~\cite{Symanzik-Varenna}), and the much more extensive results necessary to establish the existence of an interacting Euclidean field theory---even in two dimensions---remained impossible to bridge  at the time.   

Furthermore, even if one could solve that problem of constructing the probability distribution for a  Euclidean classical field,  there remained another  fundamental  question: how can one relate a classical Euclidean field to a relativistic quantum field?   In other words, as there is no Hamiltonian for the Euclidean classical field, how can one justify analytically continuing back from Euclidean time to real time?  Namely, how can one justify an inverse Wick rotation? In the rest of this work I will  call this question the \textit{reconstruction problem}; similarly I call a solution to this problem a \textit{reconstruction theorem}.

How to  prove a reconstruction theorem  baffled mathematical physicists during my student days, and it became the center of much discussion.  But there were no viable suggestions of how to resolve it. The same question also baffled particle physicists, including Schwinger\footnote{In those days,  discussions after lectures at important conferences were recorded or transcribed; then they were  published along with the lectures in the conference proceedings. In perspective, the discussions are often of greater interest than the talks. I was  grateful to learn in~\cite{Miller} about the reference~\cite{SchwingerCERN}, where one can read  this fascinating exchange. The book that I quote can be found online by following the hyperlink to CERN, and the quoted passage appears at the end of page~140.}  who regarded it as important.  This latter statement is documented in Schwinger's  response to a question asked by Pauli about understanding spin and statistics in the Euclidean framework; it followed Schwinger's talk  at a 1958 CERN conference.   ``The question of to what extent you can go backwards, remains unanswered, i.e. if one begins with an arbitrary Euclidean theory and one asks: when do you get a sensible Lorentz theory? This I do not know. The development has been in one direction only; the possibility of future progress comes from the examination of the reverse direction, and that is completely open.''  

The first big step in the reconstruction from Euclidean to real-time  was achieved in 1972 by Edward Nelson. He had been fascinated by Symanzik's Markov field approach, and formulated that framework in mathematical terms.  Nelson discovered a reconstruction theorem that starts from a Markov field and yields a scalar bosonic quantum field~\cite{NelsonReconstruction}.  He then showed that the free (Gaussian) Euclidean field satisfies the hypotheses of that construction~\cite{Nelson-Free_Field}.

%

\section{The Discovery of Reflection Positivity}
This was the situation in the summer of 1972 when I began to read the preprint of Nelson's reconstruction paper.  This led to several discussions with Konrad Osterwalder on the topic.  Not only did I wish to understand Nelson's construction, but I hoped that one could discover a more robust reconstruction method that also applied  to fermions. In the Euclidean world fermions are Grassmann, rather than abelian, so they did not appear to  fit into a Markov setting.

Furthermore, the Markov property that Nelson assumed was a strong one, and it had not been proved for the examples of two-dimensional interacting fields known at the time, although Nelson could verify it in the case of the free field. 
 Konrad was looking for an interesting project to work on, and I thought a new way to do reconstruction could be very fruitful.  I was deeply engrossed at the time in trying to understand aspects of the three-dimensional $\varphi^{4}$ theory.

Meanwhile Robert Schrader was on vacation, visiting his friend (later wife) in Germany.  
Since one could not park a car on the street for  more than one night near the Cambridge address where Robert lived, I had suggested that he leave his car in Northampton, Massachusetts.  This is a 
small college town where my wife worked, and I knew that there were no  overnight-parking laws like in Cambridge. Just after Robert returned, I planned to drive to Northampton, and from near there to fly to Chicago for a conference.  So I took Robert with me to pick up his car, and for nearly two hours we discussed the reconstruction project while I drove.  

I explained to Robert what Konrad and I had been thinking about, and I encouraged Robert to think more about the question.  But Robert immediately resisted. He too had seen Nelson's paper and thought that Nelson's construction was very natural, leaving little new to be discovered.  As the trip progressed, I was getting more and more upset---and afterward I recalled worrying about paying too little attention to the road. For not only did I feel that the question was extremely interesting, but I really hoped that Robert and Konrad could make some progress on the problem while I was in Chicago.  Finally I seemed to break through, only shortly before we arrived at Robert's car. 

About one week later, when I returned to Cambridge, Massachusetts, I found Konrad and Robert  enthusiastic, proud,  and delighted.  They thought  that they knew the key to a new reconstruction method: they had discovered the reflection positivity property!  A reflection-invariant functional $\omega$ on an algebra $\mathfrak{A}=\mathfrak{A}_{-}\otimes \mathfrak{A}_{+}$ is reflection positive on $\mathfrak{A}_{+}\ni A$ with respect to  the antilinear reflection homomorphism $\Theta$ mapping $\mathfrak{A}_{+} \mapsto \mathfrak{A}_{-}$, if $0\leqslant\omega(\Theta(A)A)$.   And in their example, the RP form provided the inner product space a Hilbert space $\mathcal{H}$ for quantum theory. Their original celebrated publication~\cite{OS1} appeared a few months later.   This is the statement for bosons; more generally one uses a twisted product $\Theta(A)\circ A$ that reduces to $\Theta(A) A$ for bosons, see~\cite{OS1,Bas-1,PAPPA}.

The RP form actually provides the quantization of the classical system.    If $\omega$ is invariant under a time-translation *-automorphism $\sigma_{t}$ for $t\in\mathbb{R}$, and $\alpha_{t}$ acts on $\mathfrak{A}_{+}$  for $t\geqslant0$, then the quantization of $\alpha_{t}$ yields a positive Hamiltonian $H$ on $\mathcal{H}$.  The importance of RP in quantum physics stems from the fact that the inner product and the Hilbert space for almost every quantum mechanics or relativistic field theory arises from the quantization of some classical system that has the  reflection-positivity property. Often one calls this property ``Osterwalder-Schrader positivity.''    

It was discovered that the original Osterwalder-Schrader paper contained an  error involving   continuity for a tensor product of distributions.  Although some persons claimed that the gap was  serious, I was sure that the problem was minor; that is what  turned out to be the case.  They corrected their method in a second paper~\cite{OS2}, where they also gave a stronger version of their reconstruction theorem--showing the equivalence of their assumptions to a modified version of Wightman's theory.  

It is also meaningful that Edward Nelson explained in his 1973 Erice Lectures~\cite{NelsonErice}, how his Markov field fits into the Osterwalder-Schrader (OS) framework, even though at that point the correction to the OS  method had not yet been found~\cite{Osterwalder}.  But by that time, Nelson had come to regard the OS reconstruction as the natural way to pass from Euclidean theory to quantum theory. 

\section{The Proliferation of Results Related to Reflection Positivity}
Since the appearance of RP in quantum field theory, it has permeated many other fields.  Here we only give an impressionistic view.   Robert Schrader left his position as postdoctoral fellow at Harvard in 1973 to go to Germany. His successor as postdoctoral fellow was  J\"urg Fr\"ohlich, who became one of the first persons to investigate RP in great detail. Right in the beginning he made an interesting study of the  generating functionals and the reconstruction theorem via RP-functional integrals, both the zero-temperature and also finite-temperature fields~\cite{JF-1}.  I refer the reader to the notes of his presentation at this meeting for several  of the many other perspectives that J\"urg pursued~\cite{FO}.

It also turns out that RP has close ties with preceding work: including with Widder's investigation of the Laplace transform in the 1930's~\cite{Widder}, as was discovered by Klein and Landau~\cite{Klein-Landau-1}. Furthermore, the reflection $\Theta$ in RP is related in various ways to the reflection $J$ in the Tomita-Takesaki theory of operator algebras~\cite{Takesaki}, as well as to other reflections in its predecessors.   Some references on these latter connections can be found in the introduction to  \cite{Bas-2} and of course the talks of Fr\"ohlich~\cite{FO} and of Longo~\cite{LO}.

\subsection{Phase Transitions and Symmetry Breaking}
Physicists believed for many years in symmetry breaking and vacuum degeneracy for certain quantum field theories.  The idea of the proof  goes back to Peierls' analysis of phase transitions in the Ising model in 1936, but it was completed mathematically by Griffiths and by Dobrushin some thirty years later.

One expects that a  $\lambda^{2}(\varphi^{2}-1/\lambda^{2})^{2}$ field theory with a ``W-shaped''  potential will also have a phase transition for the  real parameter $\lambda$ sufficently small. In this case the field will behave much like an Ising system at low temperature:  the average of the field will be approximately localized near the two minima of the potential at $\varphi=\pm\lambda$ that are separated by a large potential barrier of height $\lambda^{-2}$ at $\varphi=0$.  This became a much-sought-after result, which was announced without proof in 1973 by Dobrushin and Minlos in a widely-cited paper~\cite{DM}.

Another early application of RP that one might not have anticipated, was the important use of the Schwarz inequality arising from RP to establish global estimates of local perturbations---by reflecting the perturbation multiple times.     In particular  this was an important component  in  the first mathematical proof of the existence of symmetry breaking and vacuum degeneracy in quantum field theory.   
One used multiple-reflection RP bounds in $\lambda\varphi^{4}$ theory to estimate the deviation between local fluctuations in that model from a Peierls' type estimate, yielding a degenerate ground state for $\lambda$ sufficiently large~\cite{GJS-PT}.

\subsection{Statistical Physics}
RP played a major role in lattice statistical physics.  There is a very large literature, with early work by Fr\"ohlich, Simon, and Spencer~\cite{FSS-1} on continuous symmetry breaking, by Fr\"ohlich, Israel, Lieb, and Simon~\cite{FILS} on establishing RP, and  by Dyson, Lieb, and Simon~\cite{DLS} on establishing phase transitions in quantum spin systems. Lieb used these methods to analyze the ground state vortices in some models~\cite{Lieb}, see also \cite{Vortex}. One can consult the review of Biskup~\cite{Biskup} for much other  work.

\subsection{Relations to Mathematics}
\subsubsection{Representation Theory}
The first work relating RP to  representation theory for quantum fields arose from the desire to obtain representations of the Poincar\'e group from the analytic continuation of quantization of the representations of the Euclidean group.  After the initial work of Osterwalder and Schrader~\cite{OS1,OS2}, this was investigated abstractly by Fr\"ohlich, Osterwalder, and Seiler~\cite{FOS} and by Klein and 
Landau~\cite{Klein-Landau-1, Klein-Landau-2}.  

This was eventually developed into an entire subfield of representation theory  and stochastic analysis studied by Klein and developed extensively by Palle Jorgensen, Karl-Hermann Neeb, and  Gestur \'Olafsson, as well as their collaborators.  See~\cite{Klein-1,Palle,JNO} and the citations in these papers.

\subsubsection{Relations to PDE}
For a free scalar field, the RP property is equivalent to a statement about monotonicity of the Green's functions with respect to a change from Dirichlet  to Neumann boundary conditions.  Let $C=(-\Delta+1)^{-1}$ denote the Green's operator for the Laplacian on $\mathbb{R}^{d}$, and let $C_{D}$  denote the Green operator obtained by imposing vanishing Dirichlet boundary conditions on the time-zero hyperplane. Let  $C_{N}$ denote the corresponding Neumann Green's function for vanishing normal derivatives on the time-zero hyperplane.  Then RP is equivalent to the following statement of operator monotonicity on $L^{2}(\mathbb{R}^{d})$,
	\[
	C_{D} \leqslant C_{N}\;.
	\]
This can be seen by expressing $C_{D}$ and $C_{N}$  using the method of ``image charges,'' familiar in physics~\cite{GJ-RP}. One can generalize this to obtain a condition for RP on reflection-invariant spaces $\Sigma=\Sigma_{-}\cup\Sigma_{0}\cup\Sigma_{+}$ with an involution $\Theta$ that exchanges $\Sigma_{\pm}$ and leaves $\Sigma_{0}$ fixed. 

\subsubsection{Fourier Analysis and the Inequalities}
For the approach to classical Fourier analysis through RP, see the contribution of Frank to these proceedings~\cite{Frank}.  It is interesting that it is possible to prove Fourier bounds in subfactor theory to prove uncertainty principles, and this is related to the picture analysis that I describe in \S\ref{Sect:PictureRP} below~\cite{JLW-1,JLW-2}. This area of research on  the analytical inequalities for pictures in a non-commutative algebra  is just emerging; it seems to offer potential exciting new insights.  

\subsubsection{Relations to Tomita-Takesaki Theory and the KMS Property}
This is another enormous subject.  The reflection $\Theta$ arises in several different ways as the Tomita-Takesaki reflection operator $J$. It is  also related to the KMS property, discovered by Haag, Hugenholtz, and Winnink~\cite{HHW}.

\subsubsection{Stochastic Quantization}
The study of classical stochastic PDE's has received a great deal of attention recently.  Physicists have proposed that a classical equation with a white noise linear forcing term can be used to quantize a classical equation---the method of \textit{stochastic quantization}.  The quantum distribution arises as the limit of the distribution of solutions at infinite stochastic time.  Unfortunately reflection positivity does not hold for the distribution of the stochastic field at any finite stochastic time~\cite{Stochastic}.\footnote{We have  only proved this for a linear field. One sees in this case  that RP does hold in the infinite-stochastic-time limit. If the Wightman functions for a  non-linear  theory are continuous (as expected) in the perturbation parameter, then RP will also not hold at finite stochastic time. In that case, however, one needs an independent method to establish RP for the limit.}   Therefore it is difficult to see how to use this method to obtain a limit that satisfies RP--except when one can analyze the limit in closed form. 

\section{Some Elementary Remarks on the Mathematics of Pictures} \label{Sect:PictureRP}
 When Zhengwei Liu came to Harvard as a postdoctoral fellow in the summer of 2015, we spent a good deal of time telling each other about each of our current areas of research.  After a while  we realized that we could combine our work, by defining a pictorial framework that we called ``planar para algebras.''  This led us to a new way to think about the RP property. We implemented all this in the framework of statistical mechanics models of parafermions  on a lattice---a setting we call PAPPA.  In  these examples, we found a geometric interpretation and proof of RP~\cite{PAPPA} for a very wide class of Hamiltonians.   We then found that our models could be viewed as languages; we had  first used this principle in joint  work  with Alex Wozniakowski on quantum information~\cite{HS,PAPPA}, and later found that it provides very interesting insights into other fields of  mathematics~\cite{Picture}.  Hopefully they will also contribute in the future to physics. 
 
Parenthetically in the case of parafermions of degree $d=2$,  the fermionic case,  elementary parafermions are called Majoranas.  As quadratic functions of Majoranas represent classical and quantum lattice spins, the RP for parafermions leads to RP proofs in these cases~\cite{RP-Parafermions,Bas-1,Bas-2}.

 I believe that the simple fundamental idea embodied in our new  pictorial  approach to RP will prove fruitful in many other contexts. Let me explain some central ideas that lead to RP, without much elaboration.  
Planar algebra is based on the mathematics of a $*$-algebra of pictures.  Multiplication is vertical composition of pictures, with algebraic right-to-left order represented by top-to-bottom composition of pictures.  For a picture $T$, one represents  the adjoint map $T\mapsto T^{*}$ in the algebra by vertical reflection of the picture. If
\[
\scalebox{.7}{$
R=
\raisebox{-.9cm}
{\scalebox{1}{
\begin{tikzpicture}
\draw (0,0) --(0,1)--(1,1) --(1,0) --(0,0);
\draw (.2, -.5) -- (.2, 0);
\draw (.8, -.5) -- (.8, 0);
\draw (.2, 1) -- (.2, 1.5);
\draw (.8, 1) -- (.8, 1.5);
\node at (.5,.5) {$R$};
\end{tikzpicture}
}}\;,
\qquad\text{then}\quad
R^*=
\raisebox{-.9cm}
{\scalebox{1}{\color{blue}
\reflectbox{
\begin{tikzpicture}
\draw (0,0) --(0,1)--(1,1) --(1,0) --(0,0);
\draw (.2, -.5) -- (.2, 0);
\draw (.8, -.5) -- (.8, 0);
\draw (.2, 1) -- (.2, 1.5);
\draw (.8, 1) -- (.8, 1.5);
\node [rotate=180] at (.5,.5) {$R$};
\end{tikzpicture}
}}}\;.
$}
\]
Composed with multiplication, the vertical reflection gives an anti-linear anti-homomorphism that reverses the order of multiplication of pictures, 
\[
\scalebox{.7}{$
RK = \raisebox{-1.7cm}{
\begin{tikzpicture}
\draw (0,0) --(0,1)--(1,1) --(1,0) --(0,0);
\draw (0,1.5) --(0,2.5)--(1,2.5) --(1,1.5) --(0,1.5);
\draw (.2, 1) -- (.2, 1.5);
\draw (.8, 1) -- (.8, 1.5);
\draw (.2, 2.5) -- (.2, 3);
\draw (.8, 2.5) -- (.8, 3);
\draw (.2, -.5) -- (.2, 0);
\draw (.8, -.5) -- (.8, 0);
\node (0,0) at (.5,.5) {$R$};
\node (0,0+1.5) at (.5,.5+1.5) {$K$};
\end{tikzpicture}
}\;,
\qquad\qquad
(RK)^{*}=
\left(\raisebox{-1.7cm}{
\begin{tikzpicture}
\draw (0,0) --(0,1)--(1,1) --(1,0) --(0,0);
\draw (0,1.5) --(0,2.5)--(1,2.5) --(1,1.5) --(0,1.5);
\draw (.2, 1) -- (.2, 1.5);
\draw (.8, 1) -- (.8, 1.5);
\draw (.2, 2.5) -- (.2, 3);
\draw (.8, 2.5) -- (.8, 3);
\draw (.2, -.5) -- (.2, 0);
\draw (.8, -.5) -- (.8, 0);
\node (0,0) at (.5,.5) {$R$};
\node (0,0+1.5) at (.5,.5+1.5) {$K$};
\end{tikzpicture}
}
\right)^{*}
= \color{blue}\raisebox{-1.7cm}{
\begin{tikzpicture}
\draw (0,0) --(0,1)--(1,1) --(1,0) --(0,0);
\draw (0,1.5) --(0,2.5)--(1,2.5) --(1,1.5) --(0,1.5);
\draw (.2, 1) -- (.2, 1.5);
\draw (.8, 1) -- (.8, 1.5);
\draw (.2, 2.5) -- (.2, 3);
\draw (.8, 2.5) -- (.8, 3);
\draw (.2, -.5) -- (.2, 0);
\draw (.8, -.5) -- (.8, 0);
\node [rotate=180] at (.5,.5) {$K$};
\node [rotate=180] at (.5,.5+1.5) {$R$};
\end{tikzpicture}
}\color{black}
=K^{*}R^{*}$}\;.
\] 
It is natural to consider a second reflection $\Theta$ that is horizontal.  This reflection defines an anti-linear homomorphism of pictures,     
 \[
 \scalebox{.8}{$
\Theta(R)= \color{red} 
\raisebox{-.9cm}
{\scalebox{1}{
\begin{tikzpicture}
\draw (0,0) --(0,1)--(1,1) --(1,0) --(0,0);
\draw (.2, -.5) -- (.2, 0);
\draw (.8, -.5) -- (.8, 0);
\draw (.2, 1) -- (.2, 1.5);
\draw (.8, 1) -- (.8, 1.5);
\begin{scope}[xscale = -1];
\node  at (-.5,.5) {\reflectbox{$R$}};
\end{scope} 
\end{tikzpicture}
}}
\;,
\qquad
\color{black}
\qquad
\Theta(RK)=
\Theta\left(\raisebox{-1.7cm}{
\begin{tikzpicture}
\draw (0,0) --(0,1)--(1,1) --(1,0) --(0,0);
\draw (0,1.5) --(0,2.5)--(1,2.5) --(1,1.5) --(0,1.5);
\draw (.2, 1) -- (.2, 1.5);
\draw (.8, 1) -- (.8, 1.5);
\draw (.2, 2.5) -- (.2, 3);
\draw (.8, 2.5) -- (.8, 3);
\draw (.2, -.5) -- (.2, 0);
\draw (.8, -.5) -- (.8, 0);
\node (0,0) at (.5,.5) {$R$};
\node (0,0+1.5) at (.5,.5+1.5) {$K$};
\end{tikzpicture}
}
\right)
= \raisebox{-1.7cm}{\color{red}
\begin{tikzpicture}
\draw (0,0) --(0,1)--(1,1) --(1,0) --(0,0);
\draw (0,1.5) --(0,2.5)--(1,2.5) --(1,1.5) --(0,1.5);
\draw (.2, 1) -- (.2, 1.5);
\draw (.8, 1) -- (.8, 1.5);
\draw (.2, 2.5) -- (.2, 3);
\draw (.8, 2.5) -- (.8, 3);
\draw (.2, -.5) -- (.2, 0);
\draw (.8, -.5) -- (.8, 0);
\begin{scope}[xscale = -1];
\node at (-.5,.5) {\reflectbox{$R$}};
\end{scope}
\node at (.5,.5+1.5) {\reflectbox{$K$}};
\end{tikzpicture}
}
=\Theta(R)\Theta(K)\;.
$}
\]

The basic idea of the geometric notion of RP arises from the fact that the two $\Theta$-reflection of pictures is related to the $^{*}$-reflection of pictures by a $180$-degree rotation that we donote $\textbf{Rot}_{\pi}$.  Algebraically $\textbf{Rot}_{\pi}(\Theta(R))=R^{*}$, while in pictures,
\[
 \scalebox{.7}{$
 \textbf{Rot}_{\pi}
 \left({ \color{red}
\raisebox{-.9cm}{
{\scalebox{1}{
\begin{tikzpicture}
\draw (0,0) --(0,1)--(1,1) --(1,0) --(0,0);
\draw (.2, -.5) -- (.2, 0);
\draw (.8, -.5) -- (.8, 0);
\draw (.2, 1) -- (.2, 1.5);
\draw (.8, 1) -- (.8, 1.5);
\begin{scope}[xscale = -1];
\node  at (-.5,.5) {\reflectbox{$R$}};
\end{scope} 
\end{tikzpicture}
}}}}\right)
=
\raisebox{-.9cm}
{\scalebox{1}{
\reflectbox{\color{blue}
\begin{tikzpicture}
\draw (0,0) --(0,1)--(1,1) --(1,0) --(0,0);
\draw (.2, -.5) -- (.2, 0);
\draw (.8, -.5) -- (.8, 0);
\draw (.2, 1) -- (.2, 1.5);
\draw (.8, 1) -- (.8, 1.5);
\node [rotate=180] at (.5,.5) {$R$};
\end{tikzpicture}
}}}\;.
$}
 \]
 The horizontal reflection is the reflection for RP.  The vertical reflection, on the other hand is associated with the fact that $R^{*}R$ is positive, if the pictures have a positive expectation defining  a non-negative form that can be used to define an inner product and a Hilbert space through the GNS representation.

  In order to understand the geometric proof of RP,  one needs to analyze the horizontal multiplication of pictures. This horizontal multiplication defines the  convolution  product of operators, 
 \[
 \scalebox{.75}{
$T  \ast S =$
\raisebox{-.6cm}{\begin{tikzpicture}
\begin{scope}[shift={(16,4.5)},rotate=90,scale=.7]
\foreach \x in {0,1}{
\foreach \y in {0,1}{
\foreach \u in {0}{
\foreach \v in {1,2,3}{
\coordinate (A\u\v\x\y) at (\x+1.5*\u,\y+3*\v);
}}}}

\foreach \u in {0}{
\foreach \v in {1,2}{
\draw (A\u\v00) rectangle (A\u\v11);
\node at (0.5,6.5) {{$T$}};
\node at (0.5,3.5) {{$S$}};
}}

\draw (0,3)--++(-.5,-.5);
\draw (1,3)--++(.5,-.5);
\draw (0,7)--++(-.5,.5);
\draw (1,7)--++(.5,.5);

\draw (A0101) to [bend left=30] (A0200);

\draw (A0111) to [bend left=-30] (A0210);
\end{scope}
\end{tikzpicture}}\;.
}
\]
As there are two types of multiplication, rotation of pictures by $90$ degrees is a very important transformation.  Denote this by $\FS$; in terms of the input and output strings, this permutes them cyclically by one:
\[
\raisebox{-1.2cm}{
\scalebox{.7}[.8]{$
\raisebox{-.5cm}{
\tikz{
\node at (-1/3 + -1/6,1/2) {\size{$\mathfrak{F}_{s}$}};
\draw (-1/6,1/6) rectangle (1/6+2/3,1-1/6);
\draw (0,-1/6)--(0,1/6);
\draw (2/3,-1/6)--(2/3,1/6);
\draw (0,1+1/6)--(0,1-1/6);
\draw (2/3,1+1/6)--(2/3,1-1/6);
\node at (1/3,1/2) {\size{$T$}};
}}
=
\raisebox{-.5cm}{
\tikz{
\draw (-1/6,1/6) rectangle (1/6+2/3,1-1/6);
\draw (0,-1/6)--(0,1/6);
\draw (2/3,-1/6)--(2/3,1/6);
\draw (0,1+1/6)--(0,1-1/6);
\draw (2/3,1+1/6)--(2/3,1-1/6);
\node at (1/3,1/2) {\size{$\mathfrak{F}_{s} T$}};
}}
=
\raisebox{-.7cm}{\scalebox{.9}{
\tikz{
\draw (-1/6,1/6) rectangle (1/6+2/3,1-1/6);
\draw (0,-1/6) arc (0:-180:.25);
\draw (2/3,1+1/6) arc (180:0:.25);
\draw (-.5,-1/6)--(-.5,1+1/6);
\draw (2/3+.5,-1/6)--(2/3+.5,1+1/6);
\draw (0,-1/6)--(0,1/6);
\draw (2/3,-1/6)--(2/3,1/6);
\draw (0,1+1/6)--(0,1-1/6);
\draw (2/3,1+1/6)--(2/3,1-1/6);
\node at (1/3,1/2) {\size{$T$}};
}  }  }
$}}  \;.
\]
We named $\mathfrak{F}_{s}$  the \textit{string Fourier transform}  (SFT) as it generalizes the normal Fourier series.  
The SFT sends a product of neutral 1-qudit transformations to the convolution of two SFT's,
\[
\scalebox{.75}{$
\begin{tikzpicture}
\begin{scope}[scale=.5]
\foreach \x in {0,1}{
\foreach \y in {0,1}{
\foreach \u in {0}{
\foreach \v in {1,2,3}{
\coordinate (A\u\v\x\y) at (\x+1.5*\u,\y+3*\v);
}}}}

\foreach \u in {0}{
\foreach \v in {1,2}{
\draw (A\u\v00) rectangle (A\u\v11);
\node at (-1,5) {{$\mathfrak{F}_{s}$}};
\node at (0.5,3.5) {{$T$}};
\node at (0.5,6.5) {{$S$}};
}}

\draw (0,3)--++(-.5,-.5);
\draw (1,3)--++(.5,-.5);
\draw (0,7)--++(-.5,.5);
\draw (1,7)--++(.5,.5);

\draw (A0101) to [bend left=30] (A0200);

\draw (A0111) to [bend left=-30] (A0210);
\node at (3,5) {$=$};
\end{scope}
\end{tikzpicture}
\raisebox{.5cm}{\begin{tikzpicture}
\begin{scope}[shift={(16,4.5)},rotate=90,scale=.7]
\foreach \x in {0,1}{
\foreach \y in {0,1}{
\foreach \u in {0}{
\foreach \v in {1,2,3}{
\coordinate (A\u\v\x\y) at (\x+1.5*\u,\y+3*\v);
}}}}

\foreach \u in {0}{
\foreach \v in {1,2}{
\draw (A\u\v00) rectangle (A\u\v11);
\node at (0.5,6.5) {{$\mathfrak{F}_{s} T$}};
\node at (0.5,3.5) {{$\mathfrak{F}_{s} S$}};
};    };

\draw (0,3)--++(-.5,-.5);
\draw (1,3)--++(.5,-.5);
\draw (0,7)--++(-.5,.5);
\draw (1,7)--++(.5,.5);

\draw (A0101) to [bend left=30] (A0200);

\draw (A0111) to [bend left=-30] (A0210);
\end{scope}
\end{tikzpicture}}
\;.$}
\]
For a transformation with two input and two output strings this is a rotation by $90$ degrees.

One needs  a state on the algebra of pictures, in order to apply the GNS construction and to recover the identity of pictures with elements of a Hilbert space.  This arises from attaching the input strings to the output strings, and it defines a trace functional on the algebra of pictures.

 \subsection{The Geometric Interpretation of Reflection Positivity}
 One visualizes the new geometric proof of RP  from looking at Figure~\ref{RP-Proof}.  This illustrates how one relates the horizontal reflection $\Theta(R)$  of $R$ with its vertical reflection $R^{*}$.  The circle product has the pictorial meaning that $\Theta(R)$ and $R$ appear at the same vertical level.
\begin{figure}[h]
  \[
   \raisebox{-1.4 cm}{
 \scalebox{.5}{
 \includegraphics{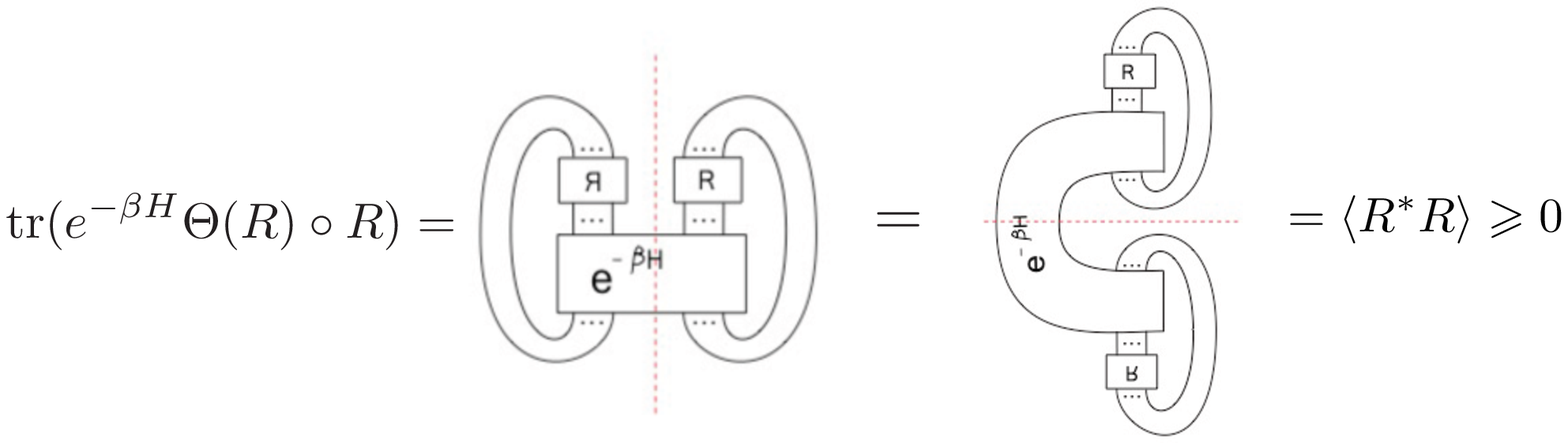}
 }}  \;.
 \]
 \caption{The Picture Proof of Reflection Positivity.}\label{RP-Proof}
 \end{figure}

  One  then rotates part of the picture to obtain an expectation of the positive Hilbert-space expectation of the quantity $R^{*}R$.   
In the picture we see that  the positive-temperature expectation of $\Theta(R)\circ R$ is given by another  expectation of $R^{*}R$, and this expectation also involves the rotation of the Hamiltonian $H$.

There are many technical points, and in~\cite{PAPPA} you will find the complete proof as well as the explanation of all the details. One must begin by proving that our  pictorial language is well-defined as a mathematical tool, and that it maps onto the objects of interest in mathematical physics. For this we implemented the parafermion algebra as an example of our model. We call the example  PAPPA (for parafermion planar para algebra). 

One must  establish the para-isotopy invariance of our pictures.  It is necessary to show that our expectation of pictures is a trace, and that it is positive. In this way one can use the GNS construction to obtain a Hilbert space representation of the pictures. 

This led us to prove a variation  of Jones' famous index theorem~\cite{Jones} for the quantization of the quantum dimension.    Finally we had to understand  how the one-string rotation $\mathfrak{F}_{s}$ of a picture  (the SFT of the picture)  is actually a generalization of the  Fourier transform of the function the picture represents. One needs to pin down the relation between  the SFT  and reflection invariance of $H$ leading to RP.  It turns out that we can  require that  $\mathfrak{F}_{s}(H)$ is positive, as it occurs in the expectation $\langle\ \cdot \ \rangle$  in~Figure~\ref{RP-Proof}.  It is not necessary that the Hamiltonian itself be positive for RP to hold, but only that its SFT is positive, This can be expressed as  $H$ having a hermitian expansion in terms of a natural basis.  
%

Having established the RP result for parafermions, one can ask how this relates to RP for classical or quantum spin systems. Many spin systems  reduce to the case for Majoranas, namely parafermions of degree $d=2$.
An extensive analysis of the relation between RP for Majoranas and  RP for classical and quantum spin systems can be found in~\cite{Bas-1}.

\section*{Acknowledgement}
I am very grateful for the hospitality and working atmosphere provided by the Research Institute for Mathematics Oberwolfach (MFO).   
This work was supported in part by the Templeton Religion Trust under Grants TRT0080 and TRT0159.

\end{talk}
\end{document}